\documentclass{article}

\usepackage{amsfonts,xr,graphicx,amsmath,amssymb}
\usepackage{srcltx}

\def\C{\mathbb{C}}
\def\S{\mathbb{S}}

\def\R{\mathbb{R}}

\def\tilde{\widetilde}
\def\epsilon{\varepsilon}

\def\n{\hfill\break} \def\al{\alpha} \def\be{\beta} \def\ga{\gamma} \def\Ga{\Gamma} \def\ro{\rho}\def\de{\delta} \def\si{\sigma}
\def\om{\omega} \def\Om{\Omega} \def\ka{\kappa} \def\la{\lambda} \def\La{\Lambda}
\def\de{\delta} \def\De{\Delta} \def\vph{\varphi} \def\vep{\varepsilon} \def\th{\theta}
\def\Th{\Theta} \def\vth{\vartheta} \def\sg{\sigma} \def\Sg{\Sigma}\def\ups{\upsilon}\def\ups{\upsilon}
\def\bendproof{$\hfill \blacksquare$} \def\wendproof{$\hfill \square$}
\def\holim{\mathop{\rm holim}} \def\span{{\rm span}} \def\mod{{\rm mod}}
\def\rank{{\rm rank}} \def\bsl{{\backslash}}
\def\il{\int\limits} \def\pt{{\partial}} \def\lra{{\longrightarrow}}
\def\noth{\varnothing}
\def\pa{\partial }
\def\ra{\rightarrow }
\def\sm{\setminus }
\def\ss{\subset }
\def\ee{\epsilon }
\def\beq{\begin{equation}}
\def\eeq{\end{equation}}
\def\ov{\over}
\def\ti{\tilde}
\def\div{\mbox{div}}
\def\supp{\mbox{supp}}
\def\tr{\mbox{tr}}
\def\const{\mbox{const}}

\numberwithin{equation}{section}
\newtheorem{theorem}{Theorem}[section]
\newtheorem{lemma}[theorem]{Lemma}
\newtheorem{proposition}[theorem]{Proposition}
\newtheorem{corollary}[theorem]{Corollary}
\newtheorem{definition}[theorem]{Definition}

\DeclareMathOperator{\curl }{curl}

\begin{document}\title{{Integral Geometry of Euler Equations }}
\author{{Nikolai Nadirashvili\thanks{Aix Marseille Universit\'e, CNRS, Centrale Marseille, I2M UMR 7373, 13453, Marseille, France, nnicolas@yandex.ru },\hskip .4 cm Serge
Vl\u adu\c t\thanks{Aix Marseille Universit\'e, CNRS, Centrale Marseille, I2M UMR
7373, 13453, Marseille, France and IITP RAS, 19 B. Karetnyi, Moscow, Russia,
serge.vladuts@univ-amu.fr} }}

\date{}
\maketitle

\def\n{\hfill\break} \def\al{\alpha} \def\be{\beta} \def\ga{\gamma} \def\Ga{\Gamma} \def\ro{\rho}\def\de{\delta} \def\si{\sigma}
\def\om{\omega} \def\Om{\Omega} \def\ka{\kappa} \def\la{\lambda} \def\La{\Lambda}
\def\de{\delta} \def\De{\Delta} \def\vph{\varphi} \def\vep{\varepsilon} \def\th{\theta}
\def\Th{\Theta} \def\vth{\vartheta} \def\sg{\sigma} \def\Sg{\Sigma}\def\ups{\upsilon}\def\ups{\upsilon}
\def\bendproof{$\hfill \blacksquare$} \def\wendproof{$\hfill \square$}
\def\holim{\mathop{\rm holim}} \def\span{{\rm span}} \def\mod{{\rm mod}}
\def\rank{{\rm rank}} \def\bsl{{\backslash}}
\def\il{\int\limits} \def\pt{{\partial}} \def\lra{{\longrightarrow}}
\def\noth{\varnothing}
\def\pa{\partial }
\def\ra{\rightarrow }
\def\sm{\setminus }
\def\ss{\subset }
\def\ee{\epsilon }
\def\beq{\begin{equation}}
\def\eeq{\end{equation}}
\def\ov{\over}
\def\ti{\tilde}
\def\div{\mbox{div}}
\def\supp{\mbox{supp}}
\def\tr{\mbox{tr}}
\def\const{\mbox{const}}

{\em Abstract.} We develop an integral geometry of stationary Euler equations
defining some function $w$ on the Grassmannian of affine lines in $\R^3$ depending on a putative compactly supported solution $(v;p)$ of the system and deduce some linear differential equations for $w$. We   conjecture that    $w=0$ everywhere and prove that this conjecture implies that  $v=0.$

\section{Introduction}
In the present paper we introduce and develop a version of integral geometry for  the steady Euler system.

More precisely, the system which we consider is as follows
\begin{equation}\label{E}
\sum_{j=1}^3\frac{\pa (v^iv^j)}{\pa x_j}+\frac{\pa p}{\pa x_i}=0\,\quad{\rm for\;\;} i=1,2,3, \eeq 
for an unknown vector field $v=v(x)=(v^1(x),v^2(x),v^3(x))$ and an unknown scalar function $p=p(x),\,x\in \R^3$; it expresses the conservation of fluid's momentum $v\otimes v +p\de_{ij}$   and  reads in a coordinate free form as
follows
\begin{equation}\label{Ed}\div (v\otimes v)+\nabla p=0\,.\eeq
Note that if we add to \eqref{E} the incompressibility condition
\begin{equation}\label{div}\div \,v=0\,,\eeq
the system \eqref{E}--\eqref{div} describes a steady state flow of the  ideal fluid.

A long-standing folklore conjecture states that a smooth compactly supported solution of  \eqref{E}--\eqref{div} should be identically zero, and this is known for Beltrami flows; see [N] and also [CC]. Let us state it explicitly:\medskip

{\bf Conjecture 1.1} {\em Let $(v;p)\in C_0(\R^3)$ be a solution of \eqref{E}--\eqref{div}. Then $v=0,\,p=0$.}\medskip

Note, however, that there do exist nontrivial Beltrami flows slowly decaying at infinity; see [EP]. Note also that nontrivial compactly supported solutions of system \eqref{E} exist, e.g., any spherically symmetric vector field $v$ is a solution of \eqref{E} for a suitable pressure $p$, but we do not know whether the   system
$$\sum_{j=1}^3v^j\frac{\pa v^i}{\pa x_j}+\frac{\pa p}{\pa x_i}=0,\;i=1,2,3 
$$ 
admits a non-trivial compactly supported solution (not satisfying \eqref{div}).\smallskip

Below we characterize the kernel of  \eqref{E} in terms of integral transforms of the quadratic forms $v^iv^j$. More precisely, given any smooth compactly supported solution $(v;p)$ of \eqref{E}, we define a smooth function $w$ on the Grassmannian  $M$ classifying lines in space using the $X$-ray transforms of $v^iv^j$ and then derive a linear differential equation for $w$. Using a Radon plane transform of $w$ we deduce one more linear homogemenous differential equation which suggests that $w=0$ everywhere. However, we are not able to deduce this fact and formulate it as a conjecture; we show that assuming the conjecture and  \eqref{div} one can deduce  Conjecture 1.1. Therefore, we put forth 

\medskip {\bf Conjecture 1.2.} {\em Let  $w$ be the function on the   Grassmannian $G_{1,3}$ of affine lines in $\R^3,$  defined below in section 3,  which depends on a compactly supported solution $(v;p)$ of \eqref{E}. Then $w=0$ everywhere.}\medskip

Note, in paricular, that this conjecture holds for any  spherically symmetric compactly supported  vector field  $v.$

\smallskip The rest of the paper is organized as follows: in Section 2 we recall some definitions and results from [S] concerning the $X$-ray transform of symmetric tensor fields. In Section 3 we define and study a smooth function $w\in  C^{\infty}(M)$   which depends on a smooth compactly supported solution $(v,p)$ of \eqref{E}. Section 4 contains a description of two invariant order 2 differential operators on $C^\infty(M)$ and a differential equation for~$w$ in terms of those operators. In Section 5 we define a plane Radon   transform for quadratic tensor fields, prove that   it  vanishes and explain why this partially confirms Conjecture 1.2. Finally, in Section 6 we deduce Conjecture 1.1 assuming Conjecture 1.2 together with  \eqref{div}.\medskip

{\em Acknowledgement}. We would like to thank A. Enciso, A. Jollivet  and V. Sharafutdinov for their aid at different stages of our work. 

\section{Tensor $X$-ray Transform }
We use throughout our paper the integral geometry of tensor fields developed in [S]
and discussed in [NSV] in its three-dimensional form. Let us give some of its points in our simple situation. For details see [S] and [NSV].

In what follows we fix a positive scalar product $\langle x,y\rangle, \, x,y\in
\R^n$. Let
$$
T{\mathbb S}^{n-1}=\{(x,\xi)\in{\mathbb R}^n\times{\mathbb R}^n:\,\|\xi\|=1,\,\langle
x,\xi\rangle =0\}\subset{\mathbb R}^n\times({\mathbb R}^n\setminus\{0\})
$$
be the tangent bundle of ${\mathbb S}^{n-1}\subset{\mathbb R}^n$.\smallskip

Given a continuous rank $h$ symmetric tensor field $f$ on ${\mathbb R}^n$, the {\it X-ray transform} of $f$ is defined for $(x ,\xi)\in{\mathbb R}^n\times ({\mathbb R}^n\setminus\{0\})$ by
\beq(If)(x,\xi)=\sum\limits_{i_1,\dots,i_h=1}^n\int\limits_ {-\infty}^\infty f_{i_1\dots i_h}(x+t\xi)\,\xi_{i_1}\dots\xi_{i_h}\,dt \label{0}
\eeq
under the assumption that $f$ decays at infinity so that the integral converges.

\smallskip
We denote by ${\mathcal S} (S^h;{\mathbb R}^n)$ the space of symmetric 
degree $h$ tensor fields with all components lying in the Schwartz space, and denote by ${\mathcal S}(T{\mathbb S}^{n-1})$ the Schwartz space on $T{\mathbb S}^{n-1}$. Below we consider only tensors from ${\mathcal S} (S^h;{\mathbb R}^n)$ and functions from ${\mathcal  S}(T{\mathbb S}^{n-1})$. For such $f\in {\mathcal S} (S^h;{\mathbb R}^n)$ we get a $C^\infty$-smooth function $\psi(x,\xi)=(If)(x,\xi)$ on ${\mathbb R}^n\times({\mathbb R}^n\setminus\{0\})$ satisfying the following conditions:
\beq\psi(x,t\xi)= \mathrm{sgn}( t) t^{h-1}\psi(x,\xi)\quad (0\neq t\in{\mathbb
R}),\quad\psi(x+t\xi,\xi)=\psi(x,\xi), \label{1}\eeq
which mean that $(If)(x,\xi)$ actually depends only on the line passing through the
point $x$ in direction $\xi$, and we parameterize the manifold of oriented lines in
${\mathbb R}^n$ by $T{\mathbb S}^{n-1}$. For $\chi(x,\xi)\in {\mathcal S}(T{\mathbb S}^{n-1})$ we can extend $\chi$ by homogeneity, setting $\chi(x,\xi)=\chi(x,\xi/\|\xi\|)$, to the open subset ${\mathcal Q}\cap \{\xi\ne 0 \}$ of the quadric
$${\mathcal Q}=\left\{(x,\xi)\in {\mathbb R}^n\times{\mathbb R}^n:\, \langle x,\xi\rangle =0\right\} \supset T{\mathbb S}^{n-1}. $$

\smallskip
Conversely, for a tensor field $f\in{\mathcal S}(S^h;{\mathbb R}^n)$, the restriction $\chi=\psi|_{T{\mathbb S}^{n-1}}$ of the function $\psi=If$ to the manifold $T{\mathbb S}^{n-1}$ belongs to ${\mathcal S}(T{\mathbb S}^{n-1} )$. Moreover, the function $\psi$ is uniquely recovered from $\chi$ by the formula
\beq\psi(x,\xi)=\|\xi\|^{h-1}\chi\Big(x-\frac{\langle x,\xi\rangle}{\|\xi\|^2} \xi,\frac{\xi}{\|\xi\|}\Big), \label{1'}
\eeq
which follows from (\ref{1}); note that $\Big(x-\frac{\langle
x,\xi\rangle}{\|\xi\|^2}\xi,\frac{\xi}{\|\xi\|}\Big)\in T\S^{n-1} \subset {\mathcal Q}$, and thus the right-hand side of \eqref{1'} is correctly defined. Therefore, the X-ray transform can be considered as a linear continuous operator $I\colon {\mathcal S}(S^h;{\mathbb R}^n)\lra{\mathcal S}(T{\mathbb S}^{n-1})$, and now we are going to describe its image and kernel.\medskip

The image of the operator $I$ is described by Theorem 2.10.1 in [S] as follows.
\smallskip 

{\bf John's Conditions.} {\em A function $\chi\in{\mathcal S}(T{\mathbb S}^{n-1})\ (n\geq3)$ belongs to the range of the operator $I$ if and only if the following two conditions hold:\smallskip

{\rm (1)} $\chi(x,-\xi)=(-1)^h\chi(x,\xi)$\textup;

{\rm (2)} The function $\psi\in C^\infty\big({\mathbb R}^n\times({\mathbb
R}^n\setminus\{0\})\big)$ defined by {\rm (\ref{1'})} satisfies the equations
\begin{equation}
\Big(\frac{\partial^2}{\partial x_{i_1}\partial\xi_{j_1}}-\frac{\partial^2}{\partial
x_{j_1}\partial\xi_{i_1}}\Big)\dots \Big(\frac{\partial^2}{\partial
x_{i_{h+1}}\partial\xi_{j_{h+1}}}-\frac{\partial^2}{\partial
x_{j_{h+1}}\partial\xi_{i_{h+1}}}\Big)\psi=0 \label{JE}
\end{equation}
written for all indices $1\leq i_1,j_1,\dots,i_{h+1},j_{h+1}\leq n$.}\medskip

Define  the symmetric inner differentiation operator $d_s=\si \nabla$ by symmetrization of the covariant differentiation operator $\nabla\colon C^{\infty} (S^h) \lra C^{\infty}(T^{h+1}),$ 
$$(\nabla u)_{i_1,\ldots  i_{h+1}} =u_{i_1\ldots i_h;i_{h+1}}=\frac{\pt u_{i_1\ldots i_h}}{\pt x_{i_{h+1}}}; $$ 
it does not depend on the choice of a coordinate system. 
\medskip 

The kernel of the operator $I$ is given by ( Theorem 2.2.1,  (1),(2) in [S]). 

\medskip {\bf Kernel of the ray transform.} {\em Let $n\ge 2$ and $h\ge 1$ be integers. For a
compactly-supported field $f\in C^{\infty}_0(S^h;\R^n)$ the following statements are equivalent:

$(1) \;If = 0;\hskip 10 cm$\smallskip

$(2)$ There exists a compactly-supported field $v \in C^{\infty}_0(S^{h-1};\R^n)$
such that its support is contained in the convex hull of the support of $f$ and
\beq\label{dv} d_sv = f. \eeq}

Note also that an inversion formula for the operator $I$ given by Theorem 2.10.2 in [S] implies that it is injective on the subspace of divergence-free (=solenoidal)
tensor fields.\medskip

{\bf The 3-dimensional case}\smallskip

For $n=3$ one notes that the tangent bundle $T\S^2$ over $\S^2$
 coincides with the homogeneous space $M=G/(\R\times{\rm SO}(2))=G'/(\R\times  {\rm O}(2))$, where $G=\R^3\rtimes{\rm SO}(3)$ is the group of proper rigid motions of $\R^3$, while $G'=G\cdot\{\pm I_3\}$ is the isometry group of $\R^3$.  Therefore  the operator $I$ for $n=3$ can be written as $I\colon{\mathcal S}(S^h;{\R}^3)\lra{\mathcal S}(M)$\smallskip

Let us define coordinates on the open subset $M_{nh}$ of $M$ consisting of
non-horizontal affine lines. Namely, $m=m\left(y_1,y_2,\alpha_1,\alpha_2 \right)$ is given by a parametric equation for a current point $A$ on $m$,
$$A=(y_1,y_2,0)+t\alpha =(y_1+\alpha_1t,y_2+\alpha_2t,t),$$
where $t$ grows in the positive direction of $m$ and thus the vector
$\alpha=(\alpha_1,\alpha_2,1)$ defines the positive direction of $m$. \medskip

We can now rewrite the above general formulas  using the
coordinates $(y_1,y_2,\alpha_1,\alpha_2)$. First we fix the following notation: \beq\label{k} k=k(\alpha_1,\alpha_2)=\sqrt{1+\alpha_1^2+ \alpha_2^2}=\sqrt{1+\|\alpha\|^2}\,;\eeq 
we will use this notation throughout the paper.\medskip

Define the diffeomorphism
$$\Phi\colon U\rightarrow{\mathbb R}^4,\quad(x,y,z,\xi)=(x,y,z,\xi_1,\xi_2,\xi_3)
\mapsto(y,\alpha)=(y_1,y_2,\alpha_1,\alpha_2),$$
on the open set $U= T{\S}^2\cap \{ \xi_3>0\}$ by
\begin{equation}y_1=x-\frac{\xi_1}{\xi_3}z,\quad y_2=y-\frac{\xi_2}{\xi_3}z, \quad\alpha_1=\frac{\xi_1}{\xi_3},\quad\alpha_2=\frac{\xi_2}{\xi_3}.\label{1.1}
\end{equation}
Then $(U,\Phi)$ is a coordinate patch on $M$; this parametrization was used by 
F. John in his seminal paper [J].\smallskip

For a function $\chi\in C^\infty(U)$, we define $\varphi\in C^\infty({\mathbb R}^4)$
by
$$
\varphi=k^{h-1}\chi\circ\Phi^{-1}.
$$
These two functions are expressed through each other by the formulas
$$
\chi(x,y,z,\xi)=\xi_3^{h-1}\varphi\left(x-\frac{\xi_1z}{\xi_3}
,y-\frac{\xi_2z}{\xi_3},\frac{\xi_1}{\xi_3},\frac{\xi_2}{\xi_3}\right),
$$
\beq\label{phi}\varphi(y,\alpha) =k^{h-1}\chi\left(y_1-\frac{\langle y,\alpha\rangle\alpha_1} {k^2},y_2- \frac{\langle y,\alpha\rangle\alpha_2}{k^2}, -\frac{\langle y,\alpha \rangle }{k}, \frac{\alpha_1 }{k},\frac{\alpha_2}{k} ,\frac{ 1}{k}\right).\eeq

If a function $\chi\in C^\infty(T{\mathbb S}^2)$ satisfies
$\chi(-x,-\xi)=(-1)^h\chi(x,\xi)$, then it is uniquely determined by
$$
\varphi=k^{h-1}\chi|_U\circ\Phi^{-1}\in C^\infty({\mathbb R}^4)\,.
$$
For a tensor field $f\in{\mathcal S}(S^h;{\mathbb R}^3)$, the function \beq
\label{phik} \varphi=k^{h-1}(If)|_U\circ\Phi^{-1}\in C^\infty({\mathbb R}^4) \eeq is
expressed through $f$ by the formula
\begin{equation}\label{1.4}
\varphi(y,\alpha)=\sum\limits_{i_1,\dots,i_h=1}^3\int\limits_{-\infty}^\infty
f_{i_1\dots i_h}(y_1+\alpha_1t,y_2+\alpha_2t,t)\,\alpha_{i_1}\dots\alpha_{i_h}dt,
\end{equation}
with $\alpha_3=1$, which easily follows from (\ref{0}).\smallskip

Let
 \beq \label{l} L=_{\rm def} {\partial^2\over \partial \alpha_2\partial y_1}- 
{\partial^2\over \partial \alpha_1 \partial y_2 } \eeq 
be the John operator. The main result of [NSV] says that for $n=3$, a function $\chi\in{\mathcal S}(T{\mathbb S}^2)$ belongs to the range of the operator $I$ for a given $h\ge 0$ if and only if the following two conditions hold:\smallskip\smallskip

{\rm (1)} $\chi(-x,-\xi)=(-1)^h\chi(x,\xi)$;\smallskip

{\rm (2)} The function $\varphi\in C^\infty({\mathbb R}^4)$ defined by {\rm
(\ref{phi})} solves the equation
\beq L^{h+1}\varphi=0. \label{1.5}\eeq
Thus $ {h^2+5h+6\over 2}$ equations \eqref{JE}  for  $n=3$  are equivalent to equation \eqref{1.5}.

\section{Function $w$}
In what follows we fix a compactly supported smooth solution $(v,p)\in
C^\infty_0(\R^3)$ of system and define a function $w\in C^\infty_0 (M)$ using the following result.
\begin{lemma}
Let $L$ be an affine plane in $\R^3$ and let $\nu_L $ be its unit normal,  then 
\beq\label{3.2} \int_{L}\langle v,z\rangle\langle v,\nu_L\rangle
\,d\sg_L =0 \eeq
for any $z\in L$ where $d\sg_L$ is the area element on $L$.
\end{lemma}

{\em Proof.} We can assume without loss of generality that $L=\{(x_1,x_2,0)\}$ and
$\nu_L =(0,0,1)=e_3,\, z=z_1e_1+z_2e_2$ for $e_1=(1,0,0),\,e_2=(0,1,0),\, e_3=(0,0,1)$. Then we have
$$ \int_{L}\langle v,z\rangle\langle v,\nu\rangle \,d\sg_L =\int (z_1v^1+z_2v^2)v^3 \,dx_1dx_2=\hskip 3 cm$$
$$=z_1\int v^1v^3 \,dx_1dx_2+z_2\int v^2v^3 \,dx_1dx_2.\hskip 1 cm$$

Note that by equation \eqref{E} with $i=1$ and  $i=2$ there holds
$$\frac{\pa (v^1v^3)}{\pa x_3} =-\frac{\pa (v^1v^2)}{\pa x_2}-\frac{\pa(v^1v^1)}{\pa x_1}-\frac{\pa p}{\pa x_1},$$
$$\frac{\pa (v^2v^3)}{\pa x_3}=- \frac{\pa (v^2v^1)}{\pa x_1}- \frac{\pa(v^2v^2)}{\pa x_2}-\frac{\pa p}{\pa x_2} ,$$
and thus we get
$$\frac{\pa }{\pa x _3}\Big(\int v^1v^3 \,dx_1dx_2\Big)=\int \frac{\pa (v^1v^3)}{\pa x_3}\,dx_1dx_2=0,
$$
$$\frac{\pa }{\pa x _3}\Big(\int v^2v^3 \,dx_1dx_2\Big)=\int \frac{\pa (v^2v^3)}{\pa x_3}\,dx_1dx_2=0.
$$
Therefore, the compactly supported functions 
$$\int v^1v^3 \,dx_1dx_2\qquad{\hbox{\rm and}}\qquad \int v^2v^3 \,dx_1dx_2$$
of $x_3$ on $\R$ are constant and thus vanish everywhere which finishes the proof.\smallskip


For any fixed value of $x_3$,  we define the vector field $v^\perp
v^3=(-v^2v^3,\,v^1v^3)$ on the plane $(x_1,x_2,x_3)$ with coordinates $\{x_1,x_2\}$ depending on $x_3$ as on a parameter, where $u^\perp=(-u_2, u_1)$ for a vector field $u=(u_1, u_2)$ on $\R^2$; note that below we use this notation for vector fields on various planes in $\R^3$. Then let us set
\beq\label{n7} F= \int_{-\infty}^\infty(-v^2v^3,\,v^1v^3)\,dx_3; \eeq 
note that $F$ is a compactly supported vector field  on the plane $\Pi_{12}=\{(x_1,x_2,0)\}$ with coordinates $\{x_1,x_2\}$ and \eqref{3.2} implies that $IF=0$. Indeed, choose $x^0 = (x^0_1, x^0_2)\in \R^2$, $0\ne \xi= ( \xi_1, \xi_2)\in \R^2$, and  let $L$ be the 2-plane through the point $x^0$ parallel to the vectors $\xi$ and $(0, 0, 1).$ Then the vector $\nu= (- \xi_1, \xi_2,0)$ is orthogonal to $L$ and the vector $\tilde \xi_a = (- \xi_1, \xi_2,a)$  is parallel to $L$ for every $a.$ By \eqref{3.2} we have
$$ \int_{L}\langle v,\tilde \xi_a \rangle\langle v,\nu\rangle
\,d\sg_L =0$$
and thus we get
$$ \int_{L}(\xi_1v^1 + \xi_2v^2 + av^3)(-\xi_2v^1 + \xi_1v^2)\,d\sg_L =0.$$
Substituting the values $a=0,\,a=1$ and taking the difference we get the equation  $IF=0$. 

Therefore, by \eqref{dv} we have
 \beq\label{3.4}d_sw_0= \nabla w_0=-F \eeq 
for a unique compactly supported smooth scalar function $w_0=w_0(x)$.\medskip

Let us fix for a moment a point $P^0=(x_1^0, x_2^0)\in \R^2$, let $r \subset\Pi_{12}= \R^2$ be a ray emanating from $P^0$ and let $e_r$ be a unit directional vector of $r$, then    in virtue of \eqref{3.4} we have
\beq\label{st}\int_{r}\langle e_{r},F\rangle \,ds_{r}=w_0(x_1^0,x_2^0)\eeq
for the  line element $ds_{r}$ of $r$. Let $H\subset \R^3$ be a half-plane perpendicular to  $\Pi_{12}$ with $\pa H=m(x_1^0, x_2^0,0,0)$, where $m(x_1^0,x_2^0,0,0)$ is the vertical line passing through the point $(x^0_1,x^0_2,0) \in \Pi_{12}$; therefore, $H$ orthogonally projects onto some ray $r$ emanating from $P_0$.
Let us consider the integral
$$ \int_{H} v^3\langle\nu_{H},v\rangle \,d\si_{H}=-\int_{r}\langle e_{r},F\rangle \,ds_{r}
$$
for the area element $d\sg_{H}$ of $H$  and a suitable unit normal $\nu_{H}$ to $H$, then  by \eqref{st}  it does not depend on $H$ for a fixed point $P^0=(x_1^0, x_2^0)$  and a fixed line $\pa H=m(x_1^0, x_2^0,0,0)$. Since the choice of a vertical line in $\R^3$ is arbitrary we see that the following definition is correct:

\begin{definition} Define
\beq\label{dw}w= -\int_{H(m)}\langle e_m,v\rangle \langle\nu_{H(m)},v\rangle \,d\si_{H(m)} \eeq
 where $H(m)$ is a half-plane with $\pa H(m)= m$ and $\nu_{H(m)}$ is the unit
normal to $H(m)$ such that the basis $\left(e_m, \nu_{m},\nu_{H(m)}\right)$ is
positively oriented for the interior unit normal $\nu_{m}$ to $m$ lying in $H(m)$.
\end{definition}

Therefore, $w$ is a compactly supported smooth function on $M$ and it can be written as $w=w(y_1,y_2,\alpha_1,\alpha_2)$ on $M_{nh}$; moreover, we get

\begin{lemma} We have
 $$w(y_1,y_2,0,0)=w_0(y_1,y_2)\,.
$$\end{lemma}

{\em Proof.} Indeed, it is sufficient to verify that 
\beq\label{P0}{\pa w\over \pa y_1}(0)=\int_{-\infty}^\infty v^2v^3 \,dx_3,\quad {\pa w\over \pa y_2}(0)=-\int_{-\infty}^\infty v^1v^3 \,dx_3, \eeq 
which is clear, since
$$ w(\de,0,0,0)=-\int_{H_{1,\de}}v^2v^3\,dx_1dx_3
=-\int_{\de}^\infty dx_1\int_{-\infty}^\infty v^2v^3 \,dx_3,$$
$$ w(0,\mu,0,0)=\int_{H_{2,\mu}}v^1v^3\,dx_2dx_3
= \int_{\mu}^\infty dx_2\int_{-\infty}^\infty v^1v^3 \,dx_3$$ for the half-planes
$$H_{1,\de}=\left\{(x_1>\de,0,x_3)\right\},\quad H_{2,\mu}=\left\{(0,x_2>\mu,x_3) \right\}.$$

Now we give two explicit formulas for $w$ which use two specific choices of~$H(m)$.
We begin by putting \beq\;k_1 =\sqrt{1+\alpha_1^2},\quad k_2
=\sqrt{1+\alpha_2^2};\eeq
recall also that $k=\sqrt{1+\alpha_1^2+\alpha_2^2}.$\medskip

Given a line $m\in M_{nh},$ let $H(m)_1$ and $H(m)_2$ be the half-planes with the border-line $m$ which are
determined by the following conditions:\smallskip

$(i)\;H(m)_1 \;{\hbox{\rm is parallel to }} x_1{\hbox{\rm -axis}},\; H(m)_2 \;
{\hbox{\rm is parallel to }} x_2{\hbox{\rm -axis;}}\hskip 3 cm $

\smallskip

$(ii)\;\langle\nu_i,e_i\rangle>0,\;i=1,2,\hskip 8 cm$\smallskip

\noindent for  $e_1=(1,0,0),\,e_2=(0,1,0)$ and 
the internal normals $ \nu_i\in H(m)_i$, $i=1,2$. \smallskip

  We have then
$$\nu_{H(m)_1}= \left(0,{ 1\over k_2},-{\alpha_2\over k_2}\right),\quad \nu_{H(m)_2}= \left(-{1\over k_1},0,{\alpha_1\over k_1} \right),
$$
and the plane $H(m)_1$ forms angle $\beta_1$ with the coordinate plane
$\Pi_{13}=\{x_2=0\}$ where $\cos\beta_1={1/k_2}$, while the plane $H(m)_2$ forms
angle $\beta_2$ with the coordinate plane $\Pi_{23}=\{x_1=0\}$, $\cos\beta_2={1/
k_1}$. Note also that we have 
$$ e_m=\frac 1k\left(\alpha_1,\alpha_2,1\right)
=\left({\alpha_1\over k},{\alpha_2\over k},{1\over k} \right)$$
for the positive unit directional vector $e_m$ of $m$.

\begin{proposition}\label{pr1}
Let $d\si_i$ be the surface area element on $H(m)_i$ and let \linebreak $l_i=y_i+x_3\alpha_i$ for $i=1,2$. Then in the introduced notation we have \beq\label{i}
\begin{aligned}[b]
\hskip-1em (i)~~ w&= -\int_{H(m)_2}\langle e_m,v \rangle\langle\nu_{H(m)_2},v\rangle
\,d\si_2 = \\ &=-\int_{-\infty}^\infty\int_{l_2 }^\infty \left(\langle e_m,v
\rangle\langle\nu_{H(m)_2},v\rangle\right)|_{(l_1,x_2,x_3)} k_1dx_2dx_3=\\& = \int_{-\infty}^\infty\int_{l_2}^\infty \frac 1k
\left((\alpha_1v^1+\alpha_2v^2+v^3)(v^1-\alpha_1v^3)\right)|_{(l_1,x_2,x_3)}dx_2dx_3
\end{aligned}
\eeq \beq\label{ii}
\begin{aligned}[b]
\hskip-5ex (ii)~~ w&= -\int_{H(m)_1}\langle e_m,v \rangle\langle\nu_{H(m)_1},v\rangle
\,d\si_1= \\ &=-\int_{-\infty}^\infty\int_{l_1}^\infty\left(\langle e_m,v
\rangle\langle\nu_{H(m)_1},v\rangle \right)|_{(x_1,l_2,x_3)} k_2\,dx_1dx_3= \\ &=
-\int_{-\infty}^\infty\int_{l_1}^\infty \frac1k\left((\alpha_1v^1+\alpha_2v^2+v^3) (v^2- \alpha_2 v^3)\right)|_{(x_1,l_2,x_3)}\,dx_1dx_3.
\end{aligned}
\eeq
\end{proposition}

{\em Proof.} This is an elementary calculation which we give only for $H(m)_1$,
since the case of $H(m)_2$ is completely similar; note only that the choice of  normals $\nu_{H(m)_1}$ and $\nu_{H(m)_2} $ comes from the orientation condition. Let us fix the values of $y_1$, $y_2$, $\alpha_1$, and $\alpha_2$, and let $H_i\supset H(m)_i$ be an affine plane containing $H(m)_i$ for $i=1,2$. Then an equation of $H_1$ is of the form
$ax_2+bx_3+c=0$, and therefore $c=-ay_2$. Since $e_m\in {\bar H_1}$ for the vector plane
${\bar H_1}$ parallel to $H_1$, we get $a \alpha_2+b=0$, and we can choose $a=1$, $b=
- \alpha_2$, so the equation takes the form
$$x_2-x_3 \alpha_2-y_2 =0,$$
and therefore $x_2=x_3\alpha_2+y_2$ on the half-plane $H(m)_1$. Since $\cos\beta_1 =
{\cos \arctan \alpha_2}= {1 \over k_2}$, we see that $d\si_{1} = k_2\,dx_1dx_3$. Then one
notes that the orthogonal projection $\pi_{13}(m)$ of $m$ on the coordinate plane
$\Pi_{13}=\{x_2=0\}$ is given by
$$\pi_{13}(m)=\Pi_{13}\cap H_2=\left\{x_1 =y_1+{\alpha_1x_3 }\right\},$$
and thus $H(m)_1$ projects onto
$$\left\{x_1>y_1+{x_3\alpha_1}\right\}\subset \Pi_{13},$$
since $\langle \nu_2,e_2\rangle>0$, which completes the proof.\medskip

The formulas \eqref{i}--\eqref{ii} are somewhat cumbersome  and use below only the following simple consequence.

\begin{corollary}\label{cor1}
In the first order of $(\alpha_1,\alpha_2)$, 
 ignoring terms with total $deg_\alpha\ge 2$, 
we have the following expressions: 
\beq\label{wI}w=\int_{-\infty}^\infty\int_{l_2}^\infty
\left(\alpha_2v^1v^2+v^1v^3+\alpha_1\left(v^1\right)^2 -\alpha_1\left(v^3\right)^2\right)|_{(l_1,x_2,x_3)}\,dx_2dx_3 ,\eeq 
\beq\label{wII}w=\int_{-\infty}^\infty\int_{l_1}^\infty\left(\alpha_2\left(v^3\right)^2-
\alpha_1v^1v^2-v^2 v^3-\alpha_2\left( v^2\right)^2\right)|_{(x_1,l_2,x_3)} dx_1dx_3 .\eeq
\end{corollary}

This corollary permits to calculate the quantities
$$\;{\pa^m w\over \pa y_1^i\pa y_2^j \pa \alpha_1^k \pa \alpha_2^l}(0),\quad i+j+k+l=m,$$
for $k+l\le 1$, and in particular, implies the following.
\begin{corollary}\label{cor2} We have
\beq\label{w12} {\pa^2w\over\pa y_2\pa\alpha_1}(0)= \int_{-\infty}^\infty\left((v^3)^2-(v^1) ^2-x_3{\pa \left(v^3v^1\right)\over\pa x_1} \right)|_{(0,0,x_3)} \,dx_3,\eeq 
\beq\label{w21} {\pa^2w\over\pa y_1\pa \alpha_2}(0)=\int_{-\infty}^\infty \left((v^2)^2 -(v^3) ^2+x_3{\pa \left(v^3v^2\right)\over\pa x_2} \right)|_{(0,0,x_3)} \,dx_3,\eeq 
\beq\label{de} {\pa^2w\over\pa y_1^2}(0)+{\pa^2w\over\pa y_2^2}(0) =\int_{-\infty}^{\infty}\left({\pa \left(v^2v^3\right)\over\pa x_1}-{\pa \left(v^1v^3\right) \over\pa x_2} \right)|_{(0,0,x_3)}\,dx_3.\eeq
 \end{corollary}\medskip

{\em Proof of \eqref{w12}.} From \eqref{wI} we have
$$w(0,y_2,\alpha_1,0)=\int_{-\infty}^\infty\int_{y_2}^\infty\left(v^1v^3+(v^1)^2
\alpha_1 -(v^3)^2\alpha_1\right)|_{(x_3 \alpha_1,x_2,x_3 )}\,dx_2dx_3, $$ 
whence we get
$${\pa w(0,0,\alpha_1,0) \over \pa y_2}=-\int_{-\infty}^\infty\left(v^1v^3+(v^1)^2 \alpha_1 -(v^3)^2\alpha_1\right)|_{(\alpha_1x_3 , 0,x_3 )}\,dx_3  $$ 
and, finally,
$${\pa^2w\over\pa y_2\pa \alpha_1}(0)=\int_{-\infty}^\infty\left((v^3)^2 -(v^1)^2  -x_3 {\pa \left(v^3v^1\right)\over\pa x_1}\right)|_{(0,0,x_3 )}\,dx_3. $$
The proof of \eqref{w21} is completely similar and that of  \eqref{de} is even simpler.

\smallskip Taking then the difference of \eqref{w12} and \eqref{w21} we get the following formula: 
\beq\label{L}{\pa^2w \over\pa y_1\pa \alpha_2}(0)-{\pa^2w \over\pa y_2\pa \alpha_1}(0) =\int_{-\infty}^\infty\big(p+(v^1)^2+(v^2)^2-(v^3) ^2
\big)|_{(0,0,x_3 )} \,dx_3 .\eeq

Indeed, we have $\;{\pa \left(v^3v^1\right)\over\pa x_1} +{\pa \left(v^3v^2\right) \over\pa x_2} =-{\pa \left(p +(v^3)^2 \right)\over\pa x_3 } \,$ by \eqref{E} and integrating
 $$\int_{-\infty}^\infty x_3 \left({\pa \left(v^3v^1\right)\over\pa x_1}  +{\pa \left(v^3v^2\right)\over\pa x_2} \right)|_{(0,0,x_3 )} dx_3=-\int_{-\infty}^\infty  {x_3 \pa \left(p +(v^3)^2 \right)\over\pa x_3 }  |_{(0,0,x_3 )} dx_3 $$
 by parts we get \eqref{L}.

\section{ Operators $P$ and $\De_M$}

Let us define first an order 2 differential operator $P$ on the space $C^{2}(M).$\smallskip

Recall that $M=G/(\R\times{\rm SO}(2))=G'/(\R\times  {\rm O}(2))$ for $G=\R^3\rtimes{\rm SO}(3)$  and $G'=G\cdot\{\pm I_3\}.$
\begin{definition}\label{d1}
Let $f\in C^{2}(M),\,m_0\in M$, and let $g(m_0)=0=(0,0,0,0)$ for $g\in G$. Then
$$Pf(m_0)=_{\rm def}Lf_g(0), $$
where $f_g(m)=f(g^{-1}(m))$ for any $m\in M$ and $L$ is defined by {\rm \eqref{l}}.
\end{definition}

\begin{lemma}\label{l1}
This definition is correct.
\end{lemma}\smallskip

{\em Proof.} We must verify that $Lf_g(0)=Lf_h(0)$ for any $g,h\in G$ such that
$g(m_0)=h(m_0)=(0)$.\smallskip

We put $ u=g^{-1}h,\,F=f_u$, and thus we have to verify that $LF(0)=LF_u(0)$ for
$u\in {\R}\times {\mathrm SO}(2)=St_0$, $St_0<G$ being the stabilizer of the vertical
line. It~is sufficient to verify the equality separately for $u\in {\R}$ and $u\in
{\mathrm SO}(2)$. It is clear for   a vertical shift $u\in {\R}$  since $L$ has constant coefficients; for a rotation $u\in {\mathrm SO}(2)$ by angle $\theta$    
 in the horizontal plane one easily calculates
$$F_u(y_1,y_2,\alpha_1,\alpha_2 )=$$
$$=F(y_1\cos \theta-y_2\sin \theta,y_2\cos \theta+y_1\sin \theta,\alpha_1\cos \theta-
\alpha_2\sin \theta,\alpha_2\cos \theta+\alpha_1\sin \theta),$$ and a simple
calculation shows the necessary equation, since we get
$${\partial^2 F_u(0)\over \partial \alpha_2\partial y_1}= \cos^2 \theta {\partial^2 F(0)\over \partial \alpha_2\partial y_1} -\sin^2 \theta{\partial^2 F(0)\over \partial \alpha_1\partial y_2}+\cos \theta \sin \theta\left({\partial^2 F(0)\over \partial \alpha_1\partial y_1} -{\partial^2 F(0)\over \partial \alpha_2\partial y_2}\right),$$
$${\partial^2 F_u(0)\over \partial \alpha_1\partial y_2}
=\cos^2 \theta {\partial^2 F(0)\over \partial \alpha_1\partial y_2}
-\sin^2 \theta {\partial^2 F(0)\over \partial \alpha_2\partial y_1}+\cos \theta \sin \theta\left({\partial^2 F(0)\over \partial \alpha_1\partial y_1} -{\partial^2 F(0)\over \partial \alpha_2\partial y_2}\right).$$ 
The proof is finished.\medskip

We can now rewrite \eqref{L} as follows
\beq\label{P} P_0w=\int_{-\infty}^\infty\big(p+(v^1)^2+(v^2)^2-(v^3) ^2 \big)\,dx_3 =\int_{-\infty}^\infty\big(p+|v|^2-2(v^3) ^2 \big)\, dx_3 \eeq
 for the operator $P_0$ being $P$ evaluated at 0, which implies that 
\beq\label{P+}Pw=\int_m\big(p+|v|^2-2\langle v,e_m\rangle ^2 \big) \,ds=\int_m\big(p+|v|^2-2v\otimes v \big)\, ds =IQ_0(m)\eeq
 for the quadratic tensor field $ Q_0=\big(p+|v|^2\big)\de_{ij}-2v\otimes v$ and any $m\in M$, since $P$ is $G$-invariant; therefore  $Pw=IQ_0$ as functions on $M$.\smallskip

We will also use the fiber-wise Laplacian $\De_M=\De_{y_1, y_2}$ acting in
tangent planes to $\S^2$; it is defined by the usual formula
$$\De_Mf(m )=\frac{\pa^2 f(m )}{\pa{y_1}^2} +\frac {\pa^2 f(m)}{\pa{y_2}^2}
$$ 
for $f\in C^2(M)$ and a vertical line $m=m(y_1,y_2,0,0)$. For any $m\in M$ the value $\De_M f(m)$ is determined by the $G$-invariance condition as for the operator $P$ above, and the rotational symmetry of $\De_{y_1, y_2}$ guarantees the correctness of that definition. Note  that the operators $P$, $\De_M$ commute and note also that \eqref{1.4} implies that for $Q\in C_0^\infty(S^h,\R^3)$ there holds a commutation rule
\beq\label{dd} I(\De Q)=\De_M(IQ).\eeq

{\em Remark 4.1.} The algebra $D_{G'}(M)$ of the $G'$-invariant differential operators on~$M$ is freely generated by $\De_M$ and $P^2$ as a commutative algebra, see [GH].\medskip

One can also give explicit formulas for $P$ and $\De_M$ in our coordinates, namely,
\beq\label{pe}P=k^2L+\alpha_1\frac{\pt}{\pt y_2}-\alpha_2\frac{\pt}{\pt y_1},\;\De_M=  k_1^2\frac{ \pt^2}{ \pt y_1^2}+k_2^2\frac{\pt^2}{\pt y_2^2} +2\alpha_1\alpha_2 {\partial^2  \over \partial y_1\partial y_2}\,.\eeq

Now we deduce the principal linear differential equation for $w$.

\begin{proposition} \label{pr4}
We have \beq\label{p4} P^2w=-4\De_M w. \eeq
\end{proposition}

{\em Proof.} We begin with the following simple result.

\begin{lemma}\label{Lk}
If $f\in C^{\infty}_0(\R^3)$ is a scalar function then $P(If)=0$.
\end{lemma}

Indeed, since $P$ is $G$-invariant, it is sufficient to verify the equation at a
single point $0\in M$ which follows from \eqref{1.5} with $h=0$.  

Lemma \ref{Lk} implies by \eqref{P+} that  
\beq\label{pwq}P^2w=PIQ_0= PI\big(\big(p+|v|^2\big)\de_{ij}-2v\otimes v\big)=-2PI(v\otimes v)\eeq  
 for a compactly supported vector field $v$ solving \eqref{E}. Moreover, we have
$$I(v\otimes v)(y_1,y_2,\alpha_1,\alpha_2)=\int_{-\infty}^{\infty}\left( {(v^3) ^2 +2v^1v^3\alpha_1+2v^2v^3\alpha_2}\right) \,{dx_3 }+O(|\alpha|^2)$$
 and thus by \eqref{de} we get
$$PI(v\otimes v)(0)=P_0I(v\otimes v)=2\int_{-\infty}^{\infty}\left({\pa \left(v^2v^3\right)\over\pa x_1} -{\pa \left(v^1v^3\right)\over\pa x_2} \right)dx_3 =2\De_{M}w(0), $$
hence $P I(v\otimes v)= 2\De_Mw $ everywhere  and $P^2w=-4\De_Mw$ by \eqref{pwq}. 

\begin{corollary}\label{cor2}The equation
\beq\label{iq0}IQ_0(m)=0,\, \forall m\in M\eeq
implies Conjecture 1.2.
 
\end{corollary}

Ideed, if $IQ_0(m)=0$ then $\De_M w(m)=-\frac14 P^2w(m)=-\frac14 PIQ_0(m)=0$ and thus $w=0$, since $w$ is  compactly supported.

\smallskip {\bf Invariant definitions and the second proof of  \eqref{p4}}

\smallskip Now let us give a description of $P$ and $\De_M$ in terms of the Lie algebra $\mathfrak g$ of $G$. We have $\mathfrak g= {\mathfrak so}(3)\oplus {\mathfrak r}(3)=\R^3\oplus\R^3$ as vector spaces, where ${\mathfrak r}(3)$ is 3-dimensional and abelian. Thus, we can write any $g\in\mathfrak g$ as $ g=(r;s) \in{\mathfrak so}(3)\oplus {\mathfrak r}(3)$, and the commutators in $\mathfrak g$ are given by
$$[(r_1;0),(r_2;0)]=(r_1\times r_2;0),\;[(0;s_1),(0;s_2)]=0,\;[(r;0),(0;s)]= (0;r \times s).
$$

Let $(R_1,R_2,R_3)$ be the standard basis of ${\mathfrak so}(3)$, and $(S_1,S_ 2, S_3)$ be that of ${\mathfrak r}(3)$. Consider the following operators on $M$:
\beq\label{0p1}\tilde \De_M = S_1^2+S_2^2+S_3^2,\; \tilde P=
 {S_1} {R_1}+  {S_2} {R_2}+  {S_3} {R_3}, \eeq
 where we denote simply  by $g$  the action on $M$ of an element $g\in U (\mathfrak g)$ of  the universal enveloping algebra  $U(\mathfrak g)$; therefore, $S_i$ acts as the infinitesimal shift in the $x_i$-direction, and $R_i$ as the infinitesimal rotation about $x_i$-axis.

\begin{proposition}\label{pr2}
We have $\tilde \De_M =\De_M$ and $\tilde P=P$.
\end{proposition}

{\em Proof.} First, the operators $\tilde\De_M$ and $ \tilde P$ are $G$-invariant.
Indeed, it follows from the rotational invariance of the quadratic form
$x_1^2+x_2^2+x_3^2$ that $\tilde \De_M$ is rotationally invariant; for translations,
the same follows from the commutation rule ${S_i}{S_j} ={S_j}{S_i}$ for
$i,j=1,2,3$.

To prove the invariance of $\tilde P$ under the $x_3$-axis rotation we verify that $\tilde P$ and ${R_3}$ commute which can be shown as follows:
$$
[ {S_1} {R_1},  {R_3}]= - {S_2} {R_1}- {S_1} {R_2},\,
[ {S_2} {R_2},  {R_3}]=  {S_1} {R_2}+ {S_2} {R_1},\,[ {S_3} {R_3}, {R_3}]= 0.
$$
Similarly we get  its  invariance  under the $x_1-$ and $x_2-$axis rotations and thus its ${\mathrm SO}(3)$-invariance, while its $x_3$-translations invariance   follows from 
$$[{S_1} {R_1}, {S_3}]= -{S_1}{S_2},\, [{S_2}{R_2},{S_3}]={S_1S_2},\, [S_3R_3, S_3]= 0.
$$
Since any line is ${\mathrm SO}(3)$-conjugate to a vertical one,  $\tilde P$ is $G$-invariant. Finally, we have $\tilde \De_M (m_0)=\De_M (m_0),\, \tilde P(m_0)=P(m_0)\,$ for $m_0=m(0,0,0,0)$, which finishes the proof, since $\De_M$ and $P$ are both $G$-invariant. Indeed, e.g., in $\tilde P(m_0)$ the terms $S_1R_1+S_2R_2$ give $L(0)$, while $S_3R_3$ vanishes since $m_0$ is invariant under both $S_3$ and $R_3.$

\medskip{\em Second proof of Proposition \ref{pr4}.} 
Let $ t\in \R$, and let $ l_t=m(t,0,0,0)$ be a vertical line in the plane
$ \Pi_{13}$; note that $\Pi_{13}=\bigcup_{t\in \R } l_t $.
\begin{lemma}\label{l2}
For $f\in C^\infty_0(M)$ we have
\beq\int_\R P^2f(l_t)\,dt =\int_\R S_2^2 R_2^2 f(l_t) \,dt.\label{l2f}
\eeq
\end{lemma}

{\em Proof.} Define the operators $A$ and $B$ by
$A= {S_1} {R_1}P,\, B= {S_3} {R_3}P;$
then
$$P^2=A+B+ {S_2}{R_2}P=A+B+{S_2R_2}\left({S_1R_1}+ S_2R_2+ S_3R_3\right).
$$ 
Since $Af(l_t)$ is a derivative of a function of $t$, while $B$ vanishes identically on $\Pi_{13}$, we get that
$$
\int_\R(Af(l_t)+Bf(l_t))dt=0.
$$
We have also 
$${S_2R_2S_1 R_1}(f(l_t))={S_1S_2R_2 R_1}(f(l_t))- {S_3S_2}{R_1}(f(l_t)),
$$
  thus the integral of the left-hand side is zero and the same is true for 
$${S_2R_2S_3R_3}(f(l_t)) ={S_3S_2R_2R_3}(f(l_t))+{S_1S_2R_3}(f(l_t));
$$
 since ${S_2R_2S_2R_2}=S_2^2R_2^2$ we get the conclusion.\medskip

\begin{lemma} \label{l3} We have
\beq\int_{\R}  R_2^2 w(l_t)\,dt =- 4\int_{\R} w(l_t)\,dt. \label{l3f}
\eeq
\end{lemma}

{\em Proof.} Let us fix a  positive constant $c<{\pi\over 2}$, and let $l^\th_t =\left({t\over\cos(\th)},0,\tan(\th),0\right)$ for any $\th$ with $|\th|<c$; therefore, $l^\th_t$ is just the line $l_t$ rotated  (in the clockwise direction) through the angle $\th$ about the origin in the plane $\Pi_{13}$, and for any $t$ we have
$$
 R_2^2w(l^\th_t)|_{\th=0}= {\pa^2\over \pa \th^2}w(l^\th_t)|_{\th=0}.
$$
Let  $e^\th_1=(\cos \th,0,-\sin\th),\;e^\th_3=(\sin\th,0,\cos \th)$ then we have
$$w(l^\th_t)={1\over\cos\th}\int_0^\infty\int_{l_t}\langle v,e^\th_1\rangle \langle v,e^\th_3\rangle |_{ \left({t\over \cos \th}+x_3 \tan\th,x_2,x_3\right)}\,dx_3 dx_2
$$
 by \eqref{i}, and if we put $t={x_1\cos \th}-x_3\sin\th$  we get that
$$\int_{\R} w(l^\th_t)\,dt={1\over \cos \th}\int_0^\infty \int_{\R^2}\langle
v,e^\th_1\rangle \langle v,e^\th_3\rangle \,dx_3 dx_2dt= \int_{x_2 >0}\langle
v,e^\th_1\rangle \langle v,e^\th_3\rangle \,dx_1dx_2 dx_3.
$$
 
Therefore,  since $\langle v,e^\th_1\rangle \langle v,e^\th_3\rangle=v^1v^3\cos
2\th+((v^1)^2-(v^3)^2){\sin 2\th\over 2}$ we have
$${\pa^2\over \pa \th^2}\int_{\R} w(l^\th_t)\,dt=\int_{ x_2 >0 }{\pa^2\over \pa \th^2} (\langle v,e^\th_1\rangle \langle v,e^\th_3\rangle) \,dx_1dx_2 dx_3= -{4}\int_{x_2 >0}
(\langle v,e^\th_1\rangle \langle v,e^\th_3\rangle)\,dx_1dx_2 dx_3
$$
\noindent and evaluating at $\th=0$ we get a proof of Lemma \ref{l3}.\medskip

 We can finish now our second proof of 
\eqref{p4}. Indeed, \eqref{l2f}--\eqref{l3f}
 imply that 
 \beq\label{4l}\int_{\R} P^2 w(l_t) dt=\int_{\R}  S_2^2 R_2^2 w(l_t) dt =-4 \int_{\R}   S_2^2w(l_t)\,dt= -4\int_{\R} \De_M w(l_t)\,dt. \eeq

If we define a function $F(x_1,x_2)$ on $\R^2$ by
$$F(x_1,x_2)= P^2w(x_1,x_2,0,0)+4\De_M w(x_1,x_2,0,0),
$$
then the integral of $F$ over the $x_1$-axis vanishes by \eqref{4l}. Changing the
coordinate system $x_1,x_2$ in $\R^2$, we get the same for the integral of $F$
over any line in the plane $\{ x_1,x_2\}$. Thus $F=0$ by Radon's theorem, and we get
the conclusion.\medskip

{\em Remark 4.2.} One can compare \eqref{p4} with results that can be deduced from \eqref{1.5} for $h=2$. A simple direct calculation using \eqref{pe} gives for $h=2$ \beq\label{pde}P^3\psi+4P\De_M \psi=0 \eeq if $\psi=IQ$ for $Q\in C_0^\infty (S^2,\R^3)$. Applying  then \eqref{pde} to  $\psi=\De_M w$ (which can be written as $\De_M w=IQ'$ for a certain $Q'$ not given here) we obtain
$P^3\De_Mw+4P\De^2_M w=\Delta_MP(P^2w+4\Delta_M w)=0  $ and thus
$P(P^2w+4\Delta_M w)=0$ which is much weaker than \eqref{p4} since the kernel of $P$ is enormous.\smallskip

However, it is possible to construct a function $u\in C_0^\infty(M)$ verifying
$$ P\De_Mu =-2\De_Mw,\quad \De_Mu=IQ_1 $$
for some $\,Q_1\in C_0^\infty(S^2,\R^3)$ and applying \eqref{pde} to $\psi=\De_Mu=IQ_1$ we get 
$$0=P^3\De_M u+4\De_M^2 P u= -2\De_M\big(P^2 w+4\De_M w\big),
$$ 
and thus we reprove \eqref{p4}. We can define $u$  similarly to \eqref{dw} as follows
$$u=\int_{H(m)}{\rm{dist}}(P,m)\left(p+\langle \nu_{H(m)},v\rangle^2\right) d\si_{H(m)},
$$
where ${\rm{dist}}(P,m)$ is the distance from a point $P\in H(m)$ to $m$.

\section{A Radon Plane Transform}
 Let us define   a Radon  tensor plane transform $J$ as follows: 
\beq \label{pt} JQ(L)= \int_L\mathrm{tr}(Q_{|L}) \,d\si_L\eeq 
for an affine plane $L\ss\R^3$ and $Q\in C^{\infty}(S^{2};\R^3)$  satisfying  
 \beq \label{est2} |Q(x)| \le C(1 +|x|)^{-2-\varepsilon},\eeq 
for some $\varepsilon>0,$  where $Q_{|L}$ is the restrictio onto $L$; definition is correct and we get a bounded linear operator 
$$J\colon {\mathcal S}(S^2;{\mathbb R}^3)\lra{\mathcal S}(\R P^3)
$$
for the manifold $\R P^3$ of affine planes $L\ss\R^3.$ 
\begin{proposition} \label{cor6}
We have $JQ_0(L)=0$.
\end{proposition}

{\em Proof.} We  have for any affine plane $L\ss\R^3$ that
$$JQ_0(L)=2\int_L\big( p+\langle v, \nu_L\rangle^2\big)\, d\sg_L=0.
 $$ 
 Indeed,  setting without loss of generality $L=\Pi_{12},\,\nu_L=e_3$  we get that
$${\pa JQ_0\over \pa x_3}(\Pi_{12})={\pa \over \pa x_3}\Big(\int_L\big( p+\langle v, \nu_L\rangle^2\big)\, d\sg_L\Big)=
$$
$$=\int_{\Pi_{12}}{\pa\big(p+( v^3)^2\big)\over\pa x_3}\, dx_1dx_2=-\int_ {\Pi _{12}}\left({\pa(v^1v^3)\over\pa x_1}+{\pa(v^2v^3)\over\pa x_2}\right)
\, dx_1dx_2=0.
$$
Therefore,  $JQ_0(\Pi_{12})$ does not depend on $x_3$ and hence equals 0. 

\medskip Let us explain in what Proposition 5.1 partially confirms   \eqref{iq0} and thus Conjecture 1.2. One can verify that the condition $JQ_0(L)=0, \forall L\in\R P^3$ is equivalent to the following  equation for the components  $\{q^{ij}\}$ of $Q_0:$
\beq\label{jq0}\sum_{ i,j=1}^3 {\pa^2 q^{ij}\over\pa x_i\pa x_j}=\De\tr\, Q_0\,,\eeq
while $IQ_0(m)=0,\, \forall m\in M$ is equivalent to the following system
$${2\,\pa q^{ij}\over\pa x_i\pa x_j}={\pa^2q^{ii}\over\pa x_j^2}+ {\pa^2 q^{jj}\over\pa x_i^2},\,1\le i<j\le 3$$
 of 3 equations and their sum gives \eqref{jq0}.
\section{A Uniqueness Theorem}
Now we can deduce    Conjecture 1.1 from Conjecture 1.2.

\begin{theorem} \label{th2}Let $(v,p)\in\C^\infty_0(\R^3)$  be a solution of \eqref{E}--\eqref {div} and let the corresponding function $w$ vanish everywhere on $M$ then $(v;p)=0$ everywhere.
\end{theorem}

{\em Proof.} For $m\in M$ denote by $t$ a vector parallel to $m$ and by $n$  a
vector perpendicular to $m,$ then the equality $w=0$ implies that
 \beq\label{int}\int_m\langle v,t\rangle\langle v,n\rangle \,ds=0. \eeq

Let $L\ss \R^3$ be an affine plane; let $ v_n=\langle\nu_L,v\rangle$ and  $v_t=v-v_n\nu_L$ be its components normal to $L$  and tangent to $L$, respectively. Then by \eqref{int} we have $IV=0$ for the vector field $V=v_n v_t$ on $L$, and hence $\curl_L V =0$ for the curl operator $\curl_L$ on $L$ which gives
 \beq\label{r} v_n\curl_L v_t - \langle v^\perp_t, \nabla_L v_n\rangle=-v_n\om_n- \langle v^\perp_t, \nabla_L v_n\rangle=0 \eeq
for the normal component $\om_n$ of $\om=\curl v$ and  the gradient operator   $\nabla_L$ on $L$. We will apply \eqref{r} to various planes $L \ss \R^3$. \smallskip

First take $L=\Pi_{12}$, then $v_t=(v^1,v^2,0),\,v_n=v^3,\,V=v^3v_t$ and we get 
\beq \label{01} v^3\curl_L v_t-\langle v^\perp_t, \nabla v^3\rangle=v^3\left({\pa v^1 \over \pa x_2} -{\pa v^2 \over \pa x_1} \right)-v^1{\pa v^3 \over \pa x_2}+v^2 {\pa v^3 \over \pa x_1} =0. \eeq 
  Let now
 $$\Om=_{\rm def}\{u\in\R^3|v(u)\ne 0,\,\om(u)\neq 0\}
$$
 and let  $D=_{\rm def}\R^3\setminus \bar \Om$.  We can suppose that  $\Om$ is not empty, since otherwise in a neighborhood of  a point $x_0$ where the maximum of $|v|$ is attained we have $\De v=-\curl\curl v+\nabla \div\, v=0$ and thus $v$ is harmonic which contradicts the maximum principle for harmonic fields. It follows then  that $v=0$ in this neighborhood and thus everywhere.

In orthonormal coordinates with $v^1(u_0)=v^2(u_0) =0$ we get for $u_0\in \Om$ that
 \beq\label{c5} v^3\left({\pa v^1 \over \pa x_2} -{\pa v^2 \over \pa x_1} \right)(u_0) =0 {\rm\;\;and\;\; therefore\;\;} \langle v,\om\rangle(u_0) =0. \eeq 
Therefore $\langle v,\om\rangle=0$ holds everywhere on $\Om$, thus on $\R^3$ and differentiating this relation in the $v$-direction  we obtain 
\beq\label{40} \langle v \nabla v,\om\rangle+ \langle v\nabla \om, v\rangle=0 ;\eeq using the commutation law
\beq\label{41} v \nabla \om= \om\nabla v \eeq 
we get then from \eqref{40} that
\beq\label{42} \langle v \nabla v,\om\rangle+ \langle\om \nabla v, v\rangle=0. \eeq
In orthonormal coordinates $x_1,x_2,x_3$ with $x_1$  directed along $v$ and ~$x_2$  directed along $\om$ at $u_0$  we can rewrite \eqref{42} as follows
\beq\label{43} {\pa v^1 \over \pa x_2}(u_0) +{\pa v^2 \over \pa x_1}(u_0)  =0, \eeq
since $v(u_0)\neq 0$ and $\om(u_0)\ne0$. Below we always use that coordinate system.  \smallskip

Moreover, since the vector $\om$ is directed along $x_2$, we get
\beq\label{45}{\pa v^1\over \pa x_2}(u_0) ={\pa v^2\over\pa x_1}(u_0) {\rm \;\; and\;\; therefore\;\;}{\pa v^1\over \pa x_2}(u_0) ={\pa v^2 \over \pa x_1}(u_0)=0. \eeq 

Let then $L=\Pi_{13}$, and thus $v_t=(v^1,0,v^3),\,v_n=v^2,\,V=v^2v_t$. Since $v_n(u_0)=0$, we get from \eqref{r} and \eqref{41} that
\beq\label{46}{\pa v^2 \over \pa x_3} (u_0)=0 {\rm \;\;and\; hence\;\; }{\pa v^3 \over \pa x_2}(u_0)  =0 {\rm \;\;and\;\; }{\pa \om^3 \over \pa x_1}(u_0)=0.\eeq
Then, differentiating \eqref{c5} with respect to $x_1$ and $x_3$ at $u_0$, we get
 \beq\label{48}{\pa \om^1 \over \pa x_1}(u_0)=0 {\rm \;\;and\;\; }
{\pa \om^1 \over \pa x_3}(u_0)=0. \eeq

Now we take $L=\{x_2+x_3=0\}$, therefore $v_n={v^2+v^3\over \sqrt 2},\nu_L=\left(0, {1\over \sqrt 2},{1\over \sqrt 2} \right)$, \linebreak $v_t=\left(v^1, {v^2-v^3\over   2},{v^3-v^2\over   2} \right)$ and  $ V={v^2+v^3\over \sqrt 2} \left(v^1, {v^2-v^3\over \sqrt 2} \right)_{\mathcal{B}}$ in the orthonormal basis ${\mathcal{B}}=\left\{e_1'=(1,0,0), e_2'= \left(0, {1\over \sqrt 2},-{1\over \sqrt 2} \right)\right\}$. Since $v_n(u_0)=0$ and the vector $v_t^\perp(u_0)$ is directed along $e'_1$, we get from \eqref{r} that
$$v^1(u_0) \left ({\pa v^2 \over \pa x_2}(u_0)+{\pa v^2\over\pa x_3}(u_0)
-{\pa v^3\over\pa x_3} (u_0)-{\pa v^3 \over \pa x_2} (u_0)\right)=0
$$
and thus
$${\pa v^2 \over \pa x_2}(u_0)+{\pa v^2\over\pa x_3} (u_0)-{\pa v^3\over\pa x_3} (0)-{\pa v^3 \over \pa x_2}(u_0)=0.$$
Therefore by \eqref{46} we get also
 \beq\label{49} {\pa v^2 \over \pa x_2}(u_0) ={\pa v^3 \over \pa x_3}(u_0).  \eeq

For $L=\{x_1+x_2=0\}$ we then have  $v_n={v^1+v^2\over\sqrt 2}$, $\nu_L=\left(  {1\over \sqrt 2},{1\over \sqrt 2},0 \right),$ \linebreak 
$v_t=\left( {v^1-v^2\over   2},{v^2-v^1\over   2},v^3\right)$ and thus
$ V={v^1+v^2\over\sqrt 2}\left({v^1-v^2\over\sqrt 2},v^3\right)_{\mathcal{B'}}$ in the plane basis 
 ${\mathcal{B'}}=\left\{\left({1\over \sqrt 2},-{1\over \sqrt 2},0 \right), (0,0,1)\right\}$. Since the vector $v_t^\perp(u_0)$ is directed along $(0,0,1)$ 
 we get from \eqref{r} and  \eqref{45}  that
$$v^1(u_0)\left(\om^2(u_0)+{\pa v^1\over\pa x_3}(u_0)+{\pa v^2\over\pa x_3} (u_0) \right)=0;
$$
therefore, we get  by \eqref{46} that
\beq\label{51} \om^2(u_0)=-{\pa v^1 \over \pa x_3}(u_0).\eeq

If a trajectory $\ga=\ga(t)$ of the flow $v$ is parametrized  by  $t,$ i.e. $ \frac{d\ga}{d t}=v,$   we have  in virtue of \eqref{41}  a differential inequality
$$|q'(t)|\leq C |q(t)|,
$$
for the function  $q(t)=_{\rm def } \om^2(\gamma(t))$ and a positive constant $C$. Therefore, if $q(0)\ne0$ then $q(t)\ne0$ for any $t\in \R$, thus any trajectory of  $v$ does not cross $\pa\Om =\pa D$ and hence stays either in $\bar\Om$ or in $\bar D$.

Using  \eqref{45}, we see that $|v|$ is constant on any trajectory $\Ga$ of the vector field $\om$ and, conversely, \eqref{46} and \eqref{48} imply that $\om$ has a constant direction on any trajectory $\ga$ of the flow $v$ and therefore  $\ga\ss \Pi_\ga$ is a plane curve for an affine plane $\Pi_\ga$. Set $\xi =\om/|\om|$, then we can define the  vector $\xi(\ga)$ for any $\ga \ss \Om$ and $\xi(\ga)\perp \Pi_\ga$. For $z\in \Om$  we get 
\beq\label{82}{\pa \xi \over\pa x_1}(z)=0\eeq
since $v(z)$ is  parallel to $x_1$, $\om(z)$ is parallel to $x_2$  and  $\xi(z)=(0,1,0)$; thus  \eqref{48} implies that
 \beq\label{80}{\pa \xi^1 \over\pa x_3}(z)=0.\eeq
 Therefore $\xi$ satisfies the Frobenius integrability condition $\langle\xi,\curl\xi \rangle =0$ and hence in a neighborhood of $z$ there exists a smooth function $U(x)$  with  $\nabla U\neq 0$ parallel to $\xi$. Moreover, $U$ is a first integral of the flow $v$ since $\pa U/\pa v=0$. Let then $S$ be a level surface of $U$ containing  $z$, then $S$ is foliated by the trajectories of $v$ and the vector field $\xi$ defines the Gauss map $\xi : S\to \S^2$. Since $\xi$ is constant on the trajectories of $v$ the image $\be=\xi(S) \ss \S^2$ is a curve or a point. Moreover,   $\be$ is orthogonal to the axis $x_1$ at $\xi (\ga)$ by \eqref{82}--\eqref{80} and   \eqref{51} implies that $\ga$ is not a straight line. Thus we can choose  $z'\in \ga, z'\ne z$ and get that  $\ga$ is orthogonal  at $\xi (\ga)$ to some line not parallel to $x_1$. Therefore $\rank\,\xi (\ga) =0,$ hence $\rank\, \xi =0$ on $S$, $\xi$ is  constant  on $S$ and thus $S$ is a plane. We see that a  neighborhood of $z$ in $\R^3$ is foliated by planes invariant under the flow $v$.

Denote by $\ga(s,t)\ss \Om$ the trajectory of $v$ passing through $z+(0,s,0)$, and let $L_s^z$ be the plane containing $\ga(s,t)$. Then  $L_s^z\perp\om (z+(0,s,0))$  and the planes $L_s^z$ are invariant under the flow $v$   and   foliate a neighborhood of $z$ in $\R^3$.  Let $\la (s,y)$ for $y\in \Pi_{13}$ be an affine function on $\Pi_{13}$ with the graph $L_s^z$ and let $l_z(y)={\pa \la\over \pa s}(0,y)$ then $l_z$ is an affine linear function. 
Denote by $G$ a connected component of $\Pi_{13}\cap \Om,\, z\in G$ then $G$ is invariant under the flow $v$. We fix some $z'\in G$ and set 
$$l=l_{z'},$$
then $l(z')=1$. 

 Let $z_1,z_2\in G$, then  some neighborhoods of $z_1$ and $z_2$ in $\R^3$ are foliated by the same set of planes  invariant under the flow $v$ and thus the sets of planes $L^{z_1}_s$ and $L^{z_2}_s$ are the same after a reparametrization. Therefore $l$  does not vanish in $G$ and we have $ l_{z_1}=l_{z_2}/l_{z_2}(z_1)$.  Set now
$$h(t)={\pa \ga(s,t)\over\pa s}{|_{s=0}}
$$
then $h(0)=(0,1,0)$ and thus by \eqref{41} the vector field $h$ is proportional to $\om$ on $\ga=\ga(0,t)$. Since $\om$ is orthogonal to $\Pi_{13}$ on $\ga$ we get that $h(y) = (0,l(y),0),$  
\beq\label{52a} {\pa v^2\over\pa x_2}(y)={\pa \ln l\over\pa v}(y) \eeq 
 for any $y\in \ga$  and $l$ does not vanish on $\ga$.

It follows by \eqref{52a} and \eqref{49} that
 \beq\label{d51}{\pa v^3\over\pa x_3}(y) ={\pa \ln l\over\pa v}(y)\eeq 
 and since $\div\, v=0$ there holds (recall that $x_1$ is directed along $v(y)$)
 \beq\label{d52}{\pa |v|\over\pa x_1}(y)=-2{\pa \ln l\over\pa v}(y) ={\pa \ln |v|\over\pa v}(y);\eeq 
hence we get that
 \beq\label{521} |v(y)|= \frac{ C_\gamma }{l^2(y)}\eeq 
along the trajectory $\ga$ for a positive constant $C_\gamma $ depending on $\ga$.  Note also that equations \eqref{52a}--\eqref{521} hold for any trajectory of $v$ in $G$ and hence by continuity in $\bar G$ outside the zero locus of $l$; in particular, we see that
\beq \label{ln0}v(y)\ne 0{ \hbox{  \rm for any }}y\in\bar G{\hbox{ \rm with }}l(y)\ne 0.\eeq
 
 Let  $z_0\in \bar G$ be a point where the function $|v|$ attains its  maximum on $\bar G$, then $z_0\in \pa G$. Indeed, if it is not the case, we have 
 $${\pa v^1 \over \pa x_3}(z_0) =0,
$$ 
 and hence $\om(z_0)=0$  by \eqref{51} which implies $z_0\in \pa G$. 

Let $\ga_1$ be the trajectory of $v$ starting from $z_1\in \pa G$ with $ v(z_1)\neq 0$, then $\om=0$ on the whole trajectory $\ga_1$ and $\ga_1\subset \pa \Om$. Therefore $\nabla b = v\times \om=0$ on $\ga_1$, where $b=p+\frac12 |v|^2 $ is the Bernoulli function (see, e.g., [AK]) and we get
\beq\label{52} \nabla p=-\frac12 \nabla |v|^2 {\rm \;\;on\;\; } \ga_1 .\eeq 
Let $y\in \ga_1$ and let  $e$ be a unit vector in $\Pi_{13}$ orthogonal to $v(y)$, then $\langle v_e(y),v(y)\rangle =0$ for $v_e= (\nabla_e  v^1,\nabla_e  v^3)$ by \eqref{51} since $\om(y)=0$. Therefore $\langle \nabla |v|^2(y),e\rangle =0$ and \eqref{52} implies that $\ga_1$ is a straight line interval $I$ which is finite  since $v$ has a compact support, $v$ vanishes at its end points and thus  $l|_I=0$ by \eqref{ln0}.

Let now  $z_0=0$ and continue to assume  that $x_1$ is directed along $v(0)\ne 0$ and $x_2$ is directed along $\om(x)\ne 0$ for some $x\in G$ (the direction of  $\om(x)$ does not depend on $x$), then $l$ is a linear function on $\Pi_{13}$ vanishing on the $x_1$-axis: $l=Cx_3$ for $ C\neq 0$. Denote now $D^+_\epsilon =B_\epsilon\cap\{  0< x_3 \}$ and $D^-_\epsilon=B_ \epsilon\cap\{0>x_3 \},$  then  we have  $D^+_\epsilon\ss G$. Indeed, first note that    $ \left(D^+_\epsilon\cap\pa G\right)\cup\left( D^-_\epsilon\cap \pa G\right)=\emptyset$  since otherwise   the trajectory $\ga_1$ starting from $z_1\in \left(D^+_\epsilon\cap \pa G\right)\cup \left(D^-_\epsilon\cap \pa G\right) $ leads to a contradiction  since $l|_{\ga_1}=0$. Moreover, for the trajectory $\al_0(t) $ of  $v^\perp =(-v^3,v^1)$ starting at 0 we have $\al_0(t)\in D^+_\epsilon$  for a small   $t>0$  and $\al_0(t)\in D^-_\epsilon$  for a small   $t<0$. Since $\om\ne 0$ on $G$ we get that $|v|$ strictly decreases along $\al_0$ by \eqref{51} ($v(0)$ being parallel to $x_1$) while $\al_0(t)$ stays in $G$ and since $|v|$ attains at $0$  its maximum in $\bar G$ we get that $D^-_\epsilon\bigcap G=\emptyset$; therefore,  $D^+_\epsilon\ss G$.

Furthermore, any trajectory $\ga_s$ of $v$ starting from the point $(0,0,s)\in D^+_ \epsilon$ with $ 0<s<\epsilon$ and a sufficiently small $\epsilon >0$  is closed. Indeed, by\eqref{521}  we may assume  that $C_{\ga_s}$ strictly increases as a function of $s\in(0,\epsilon)$ and thus  $\ga_s$ with $s\in(0,\epsilon)$ intersects the interval $(0,(0,\epsilon))$ only once. By the Poincar\'e-Bendixson theorem we get that  $\ga_s$ either

 $(i)$ tends to a limit, or

$(ii)$ tends to a limit cycle $\rho\ss G$, or 

$(iii)$ is closed. 

\noindent Since $(i)$ contradicts \eqref{521} and $(ii)$ implies that any trajectory $\ga_a$ with $s<a<\epsilon$ tends to $\rho$ which  contradicts \eqref{521} as well, we get that $(iii)$ holds. Moreover, any trajectory starting inside $D^+_\epsilon$ enters the domain$\{x_3>\de\}\cap G$ for some fixed $\de >0$
 and the union    $A=\bigcup_{0<s<\epsilon}\ga_s\ss G $ is a topological annulus.

Note now that   $C_{\ga_s}$ tends to zero for $s \to 0$  by \eqref{521}. Any trajectory $\al$ of $v^\perp$ in $A$ is orthogonal to the trajectories $\ga_s$ and thus intersects all $\ga_s$ while $|v|$ strictly decreases along $\al$ by \eqref{51} since $\om^2>0$ on $G$. If $\lambda\in (0,\epsilon)$ then  
$$\inf_{\ga_{\lambda}}|v|>|v(z)|$$
for a sufficiently small $s\in (0,\la)$ and some $z\in \gamma_s$, while the trajectory $\al_z$ of $v^\perp$ starting from $z$ intersects $\ga_{\lambda}$ and $|v|$ strictly  decreases along $\al_z$ which gives a contradiction and thus finishes the proof.

\section{Vector Analysis' Framework} 
Let us briefly discuss Conjectures 1.1. and 1.2 in terms of the  vector analysis for compactly supported tensor fields in $\R^3$. In this  section we suppose  that $v\in C_0^\infty(S^1,\R^3)$. We can rewrite  \eqref{Ed}  as follows
\beq\label{Ec}\curl(\div (v\otimes v))=0\,.\eeq
\begin{proposition} If \eqref{Ec} holds  and the corresponding function $w\in C_0^\infty(M)$   is everywhere zero on $M$  then \beq\label{sc}\Psi=\Psi(v)=_{\mathrm{def}}\si(\curl(v\otimes v))=0\,\eeq
for the symmetrization $\Psi$ of the tensor field $\curl(v\otimes v)$, i.e. $${2}\Psi^{ij}= {\epsilon_{ilm}}\frac{\pa(v^jv^l)}{\pa x_m}+ {\epsilon_{jlm}}  \frac{\pa ( v^iv^l)} {\pa x_m} \,,
$$
where $\epsilon_{ijk}$ is the standard permutation (pseudo-)tensor, giving the sign of the permutation $(ijk)$ of $(123)$ and the summation convention applies.
 \end{proposition}
{\em Proof.} For any fixed value of $x_1$,  we define the vector field 
$$Z=v^1
(v^2, v^3)=(v^1v^2,\,v^1v^3)$$ on the vertical plane ${\Pi_{23}(x_1)} =\{x_1,x_2,x_3\}$ with coordinates $\{x_2,x_3\}$ depending on $x_1$ as on a parameter.  

\smallskip We have then   $IZ(m')=w_{\nu_m}(m)=\langle \nabla w, \nu_m\rangle$ for a line $m\ss {\Pi_{23}(x_1)},$ a normal $\nu_m$ to $m$ and a line  $m'\perp m\ss {\Pi_{23}(x_1)},$  thus $IZ=0$ and hence the solenoidal component $^sZ $ equals zero, where $Z=^sZ+\,^pZ$ is the Helmholtz decomposition of the vector field $Z$. Therefore we have
$$\;\Psi^{11}=\curl Z= \curl\, ^s Z=0,$$
thus   $\Psi^{ii}=0$ for $ \,i=1,2,3$ and  rotating the coordinate system in each plane ${\{ x_i,x_j\}}$ through the angle $\frac{\pi}{4}$  we get $\Psi^{ij}=0$ for all $\,1\le i \le j\le 3.$ \medskip
 
Moreover, the proof of  Theorem 6.1  shows   that  the conditions $\Psi(v)=0$  and $\mathrm{div}\,v=0$ imply  $v=0$.\medskip

Therefore Conjectures 1.1 and 1.2   follow from \medskip

{\bf Conjecture 7.1.} {\em If $\curl(\mathrm{div}(v\otimes v))=0$ then} 
$ \si(\curl(v\otimes v))=0.$  \medskip

Another equivalent statement can be formulated as follows \medskip

{\bf Conjecture 7.2.} {\em If $I(\mathrm{div}(v\otimes v))=0$ then} 
$PI(v\otimes v)=0.$  \medskip

One can also ask whether the condition $\curl(\mathrm{div}(v\otimes v))=0$ implies that $v$ is spherically symmetric, which would grant Conjectures 1.1, 1.2, 7.1 and 7.2.

\bigskip \bigskip
\centerline{REFERENCES} 

\bigskip\noindent [AK]  V.I. Arnold, B.A. Khesin, {\it Topological Methods in Hydrodynamics}, Springer, 1998.

\medskip \noindent [CC] D. Chae, P. Constantin, {\it Remarks on a Liouville-type theorem for Beltrami flows}, Int. Math. Res. Not.  2015, 10012--10016.

\medskip\noindent [EP] A. Enciso, D. Peralta-Salas, {\it Existence of knotted vortex tubes in steady Euler flows,} Ann. Math. 175(2012), 345--367.

\medskip\noindent [GH] F. Gonzalez, S. Helgason, {\it Invariant differential operators on Grassmann manifolds,} Adv. in Math. 60(1986), 81--91.

\medskip\noindent [J] F. John, {\it The ultrahyperbolic differential equation with four independent variables,} Duke Math. J. 4(1938), no. 2, 300--322.

\medskip \noindent [N] N. Nadirashvili, {\it Liouville theorem for Beltrami flow,} Geom. Funct. An. 24 (2014), 916--921.

\medskip \noindent [NSV] N. Nadirashvili, V.A. Sharafutdinov, S. Vl\u adu\c t, {\em The John equation for tensor tomography in three-dimensions}, Inverse Problems 32(2016), 105013 (15pp), doi:10.1088/0266-5611/32/10/105013.

\medskip \noindent [S] V.A. Sharafutdinov, {\em Integral Geometry of Tensor Fields,} Utrecht: VSP, 1994.
\end{document}